\documentclass[]{article}

\input amssym.def
\def\depth{{\rm depth}\ } \def\Sat{{\rm Sat}\ 
} \def\m{{\bf{\rm m}}} \def\Hom{{\rm Hom}} \def\Ext{{\rm Ext}} 
\def\im{{\rm Im}\ } \def\Supp{{\rm Supp}\ } \def\Ass{{\rm Ass}\ } 
 \def\bp{\hbox{\it I\hskip -2pt P}} 
\def\br{\hbox{\it I\hskip -2pt R}} \def\bn{\hbox{\it I\hskip -2pt N}} 
\def\bz{\hbox{\it Z\hskip -4pt Z}} \def\bq{\hbox{\it l\hskip -5.5pt 
Q}}  \def\supp{{\rm supp \ } }

\def\ul{\underline }

\def\la{\longrightarrow}
\def\mod{\mathop{\rm mod}\nolimits}

\def\demo{\noindent{\bf Proof .-\ }}
\newtheorem{theorem}{Theorem}
\newtheorem{lemma}{Lemma}
\newtheorem{proposition}{Proposition}
\newtheorem{definition}{Definition}
\newtheorem{example}{Example}
\newtheorem{remark}{Remark}

\newtheorem{corollary}{Corollary}

\begin{document}
\begin{center}
\uppercase{{\bf On the $S_2$-fication of Some Toric Varieties }}
\end{center}
\advance\baselineskip-3pt
\vspace{2\baselineskip}
\begin{center}
{{\sc Marcel
Morales}\\
{\small Universit\'e de Grenoble I,  Institut Fourier, 
URA 188,  B.P.74, \\
38402 Saint-Martin D'H\`eres Cedex, \\
and IUFM de Lyon,  5 rue Anselme, \\ 69317 Lyon Cedex (FRANCE)}\\}
\end{center}
\vskip.5truecm\noindent

First version in 2004, revised version march 2006

{\small \sc Abstract.\ }{\small In this paper we 
prove:
\begin{enumerate}
\item {\it Some results on the Cohen-Macaulayness of the canonical 
module.} 
\item {\it We study the  
$S_2$-fication of  rings which are quotients by lattices ideals. } 
\item {\it Given a simplicial lattice ideal of codimension two $I,$
its 
Macaulayfication   is 
given explicitly from a system of generators of $I$. } 
\end{enumerate}}
\vskip1truecm\noindent
\section*{Introduction}
Let $X$ be an algebraic variety, the set of points where $X$ is not  
Cohen-Macaulay is the Non-Cohen-Macaulay locus, this set was study in
  \cite{C}. Macaulayfication is an analogous operation to resolution of singularities
 and was considered in \cite{K} where the main theorem of Macaulayfication is given.
 
 For any   affine semigroup (without torsion) $S\subset \bn^n$ let $G(S)$ be the subgroup of $\bz^n$ generated by $S$ and 
$\bar S$ be the saturation of $S$ inside $G(S)$, that is
$$\bar S=\{ m\in G(S)\ \colon\ rm\in S {\rm \ for\ some\ }r\in \bn  \},$$ it is well known that
 the normalization of the semigroup ring $K[S]$ is given by $K[\bar S]$ and Hochster proved in \cite{H} 
that $K[\bar S]$ is always a  Cohen-Macaulay ring. We have an exact sequence:
$$0\la 
K[S] \la K[\bar S] \la K[\bar S\setminus S] \la 0, $$  and $K[\bar S]$ is a Cohen-Macaulay ring
 containing $K[S]$, with the same ring of fractions. In general, the support of $K[\bar S\setminus S]$
does not coincide with the Non Cohen-Macaulay locus of $K[S]$ because  $\bar S$ is too big. Our problem consist to 
look for a ``minimal'' subsemigroup $\tilde S\subset \bar S$   containing $S$ such that $K[\tilde S ]$
 is a Cohen-Macaulay ring. 
 In  \cite{GSW} and \cite{TH} 
the authors consider 
a semigroup 
$S'\subset \bar S$ which contains $S$ such that we have an exact sequence: $$0\la 
K[S] \la K[S'] \la K[S'\setminus S] \la 0, $$
and $ \dim K[S'\setminus S] \leq n-2 .$ 
 $K[S']$  is the $S_2$-fication of $K[S]$. When $K[S']$ is a Cohen-Macaulay ring, the support of $K[\bar S\setminus S]$
 coincide with the Non Cohen-Macaulay locus of $K[S]$.  This is the case notably when 
$S$ is a simplicial semigroup.
The purpose of this paper is to give  effective methods to compute the $S_2$-fication
 for a class of toric varieties.
In the first part of this paper we consider the 
$S_2$-fication and give some general results on the Cohen 
Macaulayness of the canonical module,  one of them extends and 
improves Proposition 2.5 of \cite{G}. We also  extend and 
improve to the lattice case the above results from \cite{GSW} and \cite{TH}, 
given shorter proofs. 

In the second part we consider a codimension two simplicial toric 
ring 
$K[S]$,  and describe the Macaulayfication of this ring in terms of 
the system of generators of its ideal of definition as described   in  \cite{Mo1}, 
this ideal can be computed
 by an effective algorithm  
which works in polynomial time at very low cost. This is also implemented  in 
 my software codim2simplicial,
which computes the generators of a simplicial codimension 2 lattice ideal without using Groebner basis. 

During the meeting Current trends in Commutative Algebra held in 
 Levico,  Italy, in  
June 2002, I have submitted to Peter Schenzel, the problem developed in this paper in sections two to four,
 then we have started a joint work on this subject during more than one year. Peter Schenzel got a proof  using spectral sequences and decided to publish by himself in \cite{S3}. 
My proof developed here is completely 
different and elementary, it is a complement to Schenzel's proof.
\section {Known results on local cohomology} 
The following results are well known \cite{S1}, \cite{S2} section 
1.2.    All this 
results are also true for graded ring and modules.

Let $(R, Q)$ be a Gorenstein local ring of dimension $n$,  let $(A, \m)$ be a 
factor ring of $R$ and $M$ a finitely generated $A-$module of dimension $d$.

We recall the local duality's theorem:

 \begin{theorem} We have an isomorphism : $$H^i_{\m}(M) \simeq 
 H^i_{Q}(M)\simeq \Hom_{R}(\Ext^{ n-i }_{R}(M, R),  E(R/Q))$$
\end{theorem}
We denote by $D^i (M)$ the finitely generated $R$-module $\Ext^{ n-i 
}_{R}(M, R)$, 
and we set by $K_{M}=D^d(M)$ the canonical module.
We recall some of their 
properties:
\begin{enumerate}
        \item  For any exact sequence $$0\la M' \la M \la M''\la 0$$ we have 
a long exact 
sequence:
$$\ldots \la D^i(M'') \la D^i(M) \la D^i(M')\la D^{ i-1 }(M'') \la 
D^{ i-1 
}(M) \la D^{ i-1 }(M')\la $$

\item  $ D^i(M)=0$ for either $i>d$ or $i<0$,   $ D^d(M)$ has dimension 
$d$.  Moreover\hfill \break $\depth D^d(M)\geq \min \{ d, 2\}$,  
$D^d(M)$ satisfies
 the condition $S_{2}$  when $d\geq 2 $, and if $M$ is 
Cohen-Macaulay then so is $ D^d(M)$.

\item For all $P\in \Supp M$ we have $(D^d(M))_{P}=D^d(M_{P})$.
 \item $\dim 
D^i(M)\leq i$ for all $0\leq i <d$.  Suppose in addition that $M$ is 
equidimensional.  Then $M$ satisfies the condition $S_{k}$ if and 
only 
if $\dim D^i(M)\leq i-k$ for all $0\leq i <d$.
  \item If $M$ is 
unmixed  and  $d\geq 2 $,  then we have an exact sequence : $$0\la M \la D^d(D^d(M)) \la 
N \la 0$$ where $\dim N\leq \dim M-2$. Moreover   $M$ satisfies the 
condition $S_{2}$ if and only if  $ M $ is isomorphic to $ D^d(D^d(M))$.

\end{enumerate}

 \section {One result on the canonical module }

\begin{theorem}\label{th2} Let $(A, \m)$ be a factor ring of a Gorenstein 
local ring,  let  $M$ be a finitely generated $A-$module of dimension $ d$. 

\begin{enumerate}
\item Assume that $d\geq 3$ and  $\depth(M)>0$,  then $ \depth D^{ d-1 }(M)=0$ if and 
only if $ \depth K_{M}=2$.

\item 
 Assume 
that $d\geq 2$,  $\depth(M)= d-1$,  and  $D^{ d-1 }(M)$  has dimension $d-2$. 
 Then $ \depth D^{ d-1 }(M)= \depth K_{M}-2 $.

\end{enumerate}
In particular suppose that  $\depth(M)= d-1$,  and  $\dim D^{ d-1 }(M)=d-2$. Then  $D^{ d-1 }(M)$ is a Cohen-Macaulay module if and 
only 
  if the canonical module $K_{M}$ is Cohen-Macaulay.
\end{theorem}

Let $a\in m$ be a non zero divisor of $M$. From the exact sequence : $$ 
 0\longrightarrow M \buildrel{\times a} \over\longrightarrow M 
 \longrightarrow M/aM \longrightarrow 0$$ we  get the following long 
 exact sequence: $$ 
 0\rightarrow D^d(M) \buildrel{\times a} 
 \over\rightarrow D^d(M) \buildrel{\alpha} \over\rightarrow 
 D^{d-1}(M/aM)\buildrel{\beta} \over\rightarrow 
 D^{d-1}(M)\buildrel{\times a} \over\rightarrow D^{d-1}(M) 
 \rightarrow D^{d-2}(M/aM) \rightarrow  D^{d-2}(M) 
 \rightarrow ...
$$ From this exact sequence we get the 
short exact sequences:

\begin{equation}\label{suite1} 0\longrightarrow D^d(M) \buildrel{\times a}
\over\longrightarrow D^d(M) \longrightarrow \im \alpha \longrightarrow 
0
\end{equation}

\begin{equation}\label{suite2} 0\longrightarrow \im \alpha 
\longrightarrow D^{ d-1 }(M/aM) \longrightarrow \im \beta 
\longrightarrow 0\end{equation}
 Note that $\im 
\beta =(0:_{D^{d-1}(M)} a)$.  From the exact sequence \ref{suite2},  
we have the long local cohomology sequence: $$ 0\longrightarrow 
H^0_{\m}(\im \alpha) \longrightarrow H^0_{\m}(D^{ d-1 }(M/aM) ) 
\longrightarrow H^0_{\m}((0:_{D^{d-1}(M)} a)) \longrightarrow H^{ 1 
}_{\m}(\im \alpha) \longrightarrow 0,  $$ where $ H^0_{\m}(D^{ d-1 
}(M/aM))=H^1_{\m}(D^{ d-1 }(M/aM))=0$,  
since $\dim M/aM=d-1\geq 2$ and 
$D^{ d-1 }(M/aM) $ satisfies condition $S_{2}$,  hence the map 
$ H^0_{\m}((0:_{D^{d-1}(M)} a)) \longrightarrow H^{ 1 
}_{\m}(\im \alpha) $ is an isomorphism.

\begin{enumerate}
\item
If  $ \depth K_{M}=2$,  suppose first that $ \depth D^{ d-1 }(M)>0$,  then we can choose $a\in m$,  a non zero divisor for $ \depth D^{ d-1 }(M)$,  this will imply that $ H^{ 1 
}_{\m}(\im \alpha)=0 $ and then  $ \depth K_{M}\geq 3$. A contradiction.

If  $ \depth D^{ d-1 }(M)=0$ we have either $ \dim D^{ d-1 }(M)=0$ or not.
If $ \dim D^{ d-1 }(M)=0$ then the module $ (0:_{D^{d-1}(M)} a)$ is non null but has 
also dimension zero. If  $ \dim D^{ d-1 }(M)>0$ then choose  $a\notin \cup_{P\in 
\Ass 
 (D^{d-1}(M))\setminus \{ \m\} } P $,   we will have that $\dim (0: 
 _{D^{d-1}(M)} a)=0$ and is non null. In both cases $ H^{ 1 
}_{\m}(\im \alpha)\simeq H^0_{\m}((0:_{D^{d-1}(M)} a)) \not=0$ and so $ \depth K_{M}=2$.
\item 

  We will prove the 
  claim by induction on $d$. Remark that if $\dim M=2$,  our statement is true.  In fact following 
Section 1,  the canonical module is Cohen-Macaulay of dimension two 
and 
since by our hypothesis $D^{ 1 }(M)$ is of dimension 0,  it is Cohen 
Macaulay.
 Let $d\geq 3$, by the first claim we can assume that $\depth D^{ d-1 }(M)>0$.

 Let $a\in m$ 
be a non zero divisor for both $M$ and $D^{d-1}(M)$.
 Since $a$ is a 
non zero divisor for $D^{d-1}(M)$,  we have $\beta =0$ and we get two exact 
sequences: $$ 0\longrightarrow D^d(M) \buildrel{\times a} 
\over\longrightarrow D^d(M) \buildrel{\alpha} \over\longrightarrow 
D^{d-1}(M/aM)\la 0$$ $$0\longrightarrow D^{d-1}(M)\buildrel{\times a} 
\over\longrightarrow D^{d-1}(M) \longrightarrow D^{d-2}(M/aM) 
\longrightarrow 0
$$
It then follows that $M/aM$ satisfies the induction hypothesis.  
Hence  
$\depth D^{d-1}(M/aM)=\depth  D^{d-2}(M/aM)+2 $,  the above two short exact sequences imply then that $\depth D^{d}(M)=\depth  D^{d-1}(M)+2 $.

\end{enumerate}
This ends 
the proof of the theorem.
As a consequence of the proof we have:
\begin{corollary} Let $(A, \m)$ be a factor ring of a Gorenstein 
local ring,  let $M$ be a finitely generated $A-$module with $\dim M= 
d\geq 3$ and  $\depth(M)>0 $.
If $\dim D^{d-1}(M)> 0$,  let $a\in m$ be a non zero divisor of $M$ and 
$a\notin \cup_{P\in \Ass (D^{d-1}(M))\setminus \{ \m\} } P $, 
then $K_{M}/aK_M$ is isomorphic to $K_{M/aM}$ if and only if $\depth K_M\geq 3$. In particular
 if 
$K_{M}$ is a Cohen-Macaulay module then $K_{M/aM}$ is a Cohen-Macaulay 
module. 
\end{corollary}
\demo With the above notations,  $K_{M}/aK_M$ is isomorphic to $K_{M/aM}$ if and only 
if $\im \beta=(0:_{D^{d-1}(M)}a)=0$. By our choice of $a$,  we have that 
$\dim (0:_{D^{d-1}(M)}a)=0$,  so $H^0_{\m}((0:_{D^{d-1}(M)}~a))=(0:_{D^{d-1}(M)}a)$ 
and then  $\im \beta=0$ if and only if 
 $H^{ 1 }_{\m}(\im 
\alpha) =0$,  this is equivalent to $\depth K_M\geq 3$.

\begin{example} (See also\cite{TH}) Consider the semigroup in $S\subset \bn^3$ generated 
by the elements $(3, 0, 0), (2, 1, 0), $ $(0, 3, 0), (3, 0, 1), (2, 1, 1), (0, 3, 1)$,  the 
semigroup ring $K[ S] $ has dimension 3,  codimension 3,  and   $\depth K[ S]=2$,  the ring $K[ S]$ 
satisfies the condition $S_{2}$ of Serre so it is isomorphic to 
$D^{3}(D^{3}(K[ S]))$. The canonical module $ D^{3}(K[ S]) $ is not 
Cohen-Macaulay.  Remark that we have $ \dim D^{2}(K[ S])= 0 .$
        
\end{example}

 \section { $S_2$-fication of unmixed modules}
        Let $(A, \m)$ be a noetherian local ring,  (resp.  graded),  
 quotient of a Gorenstein local ring (resp.  
graded Gorenstein ring) and $M$ be an $A-$module of dimension $d$.

We recall that if $M$ is unmixed,  the module  $D^d(D^d(M))$  
satisfies the condition $S_{2}$ and we have an exact sequence :
 
$$0\longrightarrow {M }\longrightarrow {D^d(D^d(M)) }\longrightarrow 
{M''
}\longrightarrow 0$$ with $\dim M''\leq d-2$.
Moreover if there exist an $A-$module $M'$ of dimension $d$,  satisfying the 
condition $S_{2}$ and an exact sequence : $$0\longrightarrow {M 
}\longrightarrow {M' }\longrightarrow {M'/M }\longrightarrow 0$$ with 
$\dim M'/M\leq d-2$,  then $M'\simeq D^d(D^d(M))$. The $A-$module $M'$ is the $S_2$-fication of $M$
 and if $M'$ is a Cohen-Macaulay module it is a Macaulayfication of $M$.

 \begin{lemma}Set $M':=D^d(D^d(M))$. Assume 
that $M$ is unmixed not 
 satisfying the condition $S_{2}$,  then:
 \begin{description}
 \item[A)] ¥The canonical module $K_{M}= D^d(M)$ is a Cohen-Macaulay 
 module if and only if $M'$ it is.
 
 \item[B)] If $K_{M}$ is a Cohen-Macaulay module,  then:

 $H^{ i-1 }_{\m}¥(M'/M)\simeq H^i_{\m}¥(M)$ for $i=1, \ldots,  
d-1$. 
In particular $\depth (M'/M)=\depth M-1$ and $\dim M'/M= \max \{ i\leq 
d-2\  / \ 
H^{ i+1 }_{\m}¥(M)\not= 0\} $

As a special case $M'/M$ is a Cohen-Macaulay Module if and only if  
only one 
of the local cohomology modules $ 
H^{ i }_{\m}¥(M),  i\leq n-1,  $  does not vanish.

In this case the Matlis dual $D^i(M)$ of $H^{ i }_{m}¥(M)$ is a Cohen 
Macaulay module of dimension $i-1$.

 In particular if $\depth M= d-1, $ then $M'/M$ is a 
Cohen-Macaulay Module of dimension $d-2$,  and $D^{ d-1 }(M)$  is a 
Cohen-Macaulay module
 of dimension $d-2$.
 
\item[C)] The Non-Cohen-Macaulay locus of $M$ 
is given by $\Supp (M'/M)$.
 
 \end{description}¥

  \end{lemma}
  \demo 
   \begin{description}
        \item[A)]  Since $\dim M'/M\leq n-2$ we have $D^d(M)\simeq D^d(M')$,  if $M'$ 
  is Cohen-Macaulay,  then $ D^d(M')$ is a Cohen-Macaulay module,  hence 
  the canonical module $K_{M}$ is Cohen-Macaulay of dimension $d$.  The 
  converse follows since $M'\simeq D^d(D^d(M'))$
   
        \item[B)]  From the long exact sequence of the local cohomology associated 
to the sequence:
 $$0\longrightarrow {M 
}\longrightarrow {M' }\longrightarrow {M'/M }\longrightarrow 0$$ 
 we get $H^{ i-1 }_{\m}¥
(M'/M)\simeq 
H^i_{\m}¥(M)$ for $i=1, \ldots,  d-1$,  which implies B).

 \item[C)] The above exact sequence is still exact by localization on 
any 
 prime ideal $P$; on the other hand $K(K(M))_{P}=K(K(M_{P}))$ and 
 recall that if $M$ is Cohen-Macaulay then the natural map 
 $M\longrightarrow K(K(M))$ is an isomorphism.  It 
 follows that the Non-Cohen-Macaulay locus of $M$ is given by $\Supp 
(M'/M)$.

   \end{description}
 When $M$ is unmixed we get a new version of theorem \ref{th2}:  
 \begin{theorem}\label{th3} Let $M$ be unmixed of dimension $d$,  not 
 satisfying the condition $S_{2}$,  and $\depth M=d-1$,  then $\dim D^{ 
 d-1 }(M )=d-2$,  and   $\depth D^{ d-1 }(M )=\depth K_{M}-2$. In particular  $D^{ d-1 }(M )$ is Cohen-Macaulay 
  if and only if 
 $K_{M}$ is Cohen-Macaulay. 
\end{theorem}
\demo In regard of Theorem 2,  we need only to prove that $\dim D^{ d-1 }(M 
)=d-2$.

Set $M'= D^d(D^d(M))$,  from the exact sequence : $$0\longrightarrow {M 
}\longrightarrow {M' }\longrightarrow {M'/M }\longrightarrow 0$$ with 
$\dim M'/M\leq d-2$,  we have $D^d(M)\simeq D^d(M')$.  Since $M'$ 
satisfies $S_{2}$ we have $\dim D^{ d-j }(M')\leq d-j-2$,  for all 
$0\leq j \leq d-1$.

 Assume that $\dim M'/M=d-k< d-2$,  
for some $3\leq k \leq d-1$,  then from the long exact sequence 
associated to the above short exact sequence we have: $$0\la D^{ d-k 
}(M'/M)\la D^{ d-k }(M')\la 0, $$ this carries a contradiction since 
$d-k=\dim M'/M= \dim D^{ d-k }(M')=d-k\leq d-k-2$.

  From the exact 
sequence $$0\la D^{ d-1 }(M') \la D^{ d-1 }(M) \la D^{ d-2 }(M'/M)\la 
D^{ d-2 }(M')\la 0$$ we get $\dim D^{ d-1 }(M )= \dim D^{ d-2 }(M'/M)= 
d-2$ since $\dim D^{ d-j }(M')\leq d-j-2$,  for all $0\leq j \leq d-1$.

\begin{remark}\label{rem-syz}  Let $M$ be a finitely generated graded 
module over a ring of polynomials,  with $\dim M=d,$ and $  \depth M=d-1$.  It 
is well known that if $0\la G \buildrel{ \phi }\over\la F \la\ldots$ 
 is the last term of the minimal syzygies of $M$ then $\ldots\la F \buildrel{ 
\sigma }\over\la G \la D^{ d-1 }(M) \la 0$ is a presentation of $D^{ d-1 }(M) $,  
where $\sigma$ is the matrix transpose of $\phi$.
\end{remark}¥

\begin{example}
Let $A$  be the affine ring of the projective surface in $P^4$ defined 
parametrically by: $$a=s^4+t^4,  b=s^2tu,   c=s^3t,  d= st^3 ,  e=su^3$$ then
 $depth A=2$ and $\sigma= 
(a, c^2+d^2, ce,  b^3)$.  It follows that $D^{ 2 }(A)$ is Cohen-Macaulay of dimension 1, 
and the 
$S_{2}-$fication is in fact a Macaulayfication.  It is not difficult 
to check that $$A'=K[s^4+t^4,  s^2tu,  s^3t,  st^3 , su^3, s^2t^2]$$ is the  
Macaulayfication of  $A$.
\end{example}
\begin{example}
Let $A$ be the affine ring of the projective surface in $P^4$ defined 
parametrically by: $$a=s^4,  b=s^3t+u^4,  c=s^2t^2,  d= su^3 ,  
e=t^2u^2$$ a quick 
computation with Macaulay, if $char(K)\not=2,3$, gives that $depth A=2$ and $\sigma$ is 
given by $$ 
\pmatrix{0 &-ae &d^2 &-c &b &0 &0 \cr 2ab^3-2a^2bc+d^4 
&3/2ab^3+1/2a^2bc 
&a^2b^2-a^3c &-3/2bd^2-a^2e &-1/2ad^2 &e &6c }, $$ in this case $\m$ 
is an 
associated prime ideal of the ideal generated by the entries of the 
second row of $\sigma$,  but again using Macaulay we get that that 
$D^{ 
2 }(A)$ is Cohen-Macaulay  of dimension 1.  Also in this example the 
$S_{2}-$fication is in fact a Macaulayfication. We can 
check that $$A'=K[s^4,  s^3t+u^4,  s^2t^2,   su^3 ,  t^2u^2, s^2u^2]$$ is the  
Macaulayfication of  $A$.
\end{example}

\section{Lattice and toric ideals}

Let $R=K[ x_{1},  \ldots ,  x_{m}] $ be a polynomial ring,  $L\subset 
\bz 
^m$ a lattice of rank $r$.  We assume that $L$ is a positive lattice,  
that is,  every non zero vector in $L$ has positive and negative coordinates.
 We 
can write every vector ${\bf u}$ in $\bz ^m$ uniquely as ${\bf 
u}={\bf u}_+-{\bf u}_-$,  where ${\bf u}_+$ and ${\bf u}_-$ are 
non-negative and have disjoint support.  Set $I_{L}\subset R$ be the 
ideal generated by all the binomials ${\underline \bf x }^{\bf u_+}- 
{\underline \bf x }^{\bf u_-}$,  where ${\bf u}$ runs over all vectors 
of $L$.  $I_{L}\subset R$ is called a lattice ideal associated to $L$.  
Let $\Sat L= \{ {\bf u}\in \bz^m \mid k {\bf u} \in L { \ \rm for\ 
some\ } k\in \bz \}$.  The group $\bz ^m/L$ is a finitely generated 
abelian group,  and $I_{L}\subset R$ is a prime ideal if and only if 
$\bz ^m/L$ has no torsion. 
 We quote the following theorem from the proof of 
Corollaries 2.2 and 2.5 of \cite{ES}:

\begin{theorem} Let $K$ be an algebraically closed field of 
any characteristic $p\geq 0$.  The ideal $I_{L}\subset R$ is always 
unmixed.  Moreover any $x_{i}$ is a non zero divisor modulo $I_{L}$.
\end {theorem}

When the ideal $I_{L}\subset R$ is prime it 
is called toric.  In the toric case the lattice $L$ is usually  viewed as the 
lattice of the relations of a finitely generated semigroup $S\subset 
\bn^n$.  In  general   we have an isomorphism 
$\bz^m/L \longrightarrow \bz^d\oplus H$,  where $H$ is a finite group,  
the images $ { \bf a_{1} }\ldots, { \bf a_{m} }$ of the canonical basis 
of $\bz^m$ under this isomorphism generate a finitely generated 
semigroup $S\subset \bz^n\oplus H$,  which generates $G(S):=\bz^d\oplus H$.
 In fact $K[ S] :=R/I_{L}= 
\oplus_{g\in S}¥ K\underline t^{ g_{1}¥ }\underline u^{ g_{2}¥ }$,  
where $g_{1}\in \bz^n,  g_{2} \in H$ and $g=(g_{1}, g_{2})$.

We set $\tilde S$ the projection of $ S$ in $\bz^d$,  let ${ \cal C 
}_{S}$ be the cone generated by $\tilde S$ in $\bq^d$,  and $F_{1}¥,  
F_{2}, \ldots,  F_{l}$ its faces of dimension $d-1$.  Let $S_{i}= \{ 
x-y; x,  y\in S,  \tilde y\in \tilde S\cap F_{i}\}\subset \bz^n\oplus 
H$.

Let $L$ be any  lattice,  corresponding to the semigroup 
$S\subset G(S)$, as in \cite{GSW},  we will define another semigroup 
$S'=\cap S_{i}¥\subset G(S)$ such that the semigroup ring $K[ S'] $ is 
the $S_{2}$-fication of the semigroup ring $K[ S] :=R/I_{L}$.  
Moreover if $S$ is simplicial then $K[ S'] $ is the Macaulayfication 
of the semigroup ring $K[ S] :=R/I_{L}$. This extends to the lattice 
case a theorem of \cite{GSW}.  First we extends some preliminary results from 
\cite{TH},  to the lattice case,  the proofs are very similar and we let it 
to the reader.

\begin{enumerate}
\item Let $\bar S= \{ z\in G(S) ,  \exists p\in \bn^*,  pz\in S\}$,  then 
$K[ \bar S] $ is the normalization of $K[ S] $. 
 \item A semigroup 
$S\subset \bz^l\oplus H$, is called 
standard if the following conditions are satisfied:
\begin{enumerate}
\item $\bar S=G(S)\cap \bn^l\oplus H, $

\item $S_{(i)}\not= S_{(j)}$ for $i\not= j$,  where $S_{(i)}=\{ x\in S; 
x_{i}=0\}$ and $x=(x_{1}, \ldots, x_{l}, h)$,  with $h\in H$.

\item ${\rm rank}_{\bz}¥ G(S_{(i)})={\rm rank}_{\bz}¥ G(S) -1,\   i=1, \ldots, l.$
\end{enumerate}¥
 
By the Hochster's transformation,  see \cite{TH},  there is an standard 
semigroup $T(S)$ isomorphic to $S$,  also by this transformation 
$T(S')=T(S)'$.  So we can assume that our semigroup $S$ is standard.

\item 
The 
polynomial ring $R$ has two gradings,  it is $G(S)=\bz^m/L=\bz^d\oplus 
H$-graded: two monomials ${\underline x}^\alpha,  {\underline x}^\beta 
$ have the same grading if and only if the vector $\alpha -\beta \in 
L$,  the lattice ideal $I_{L}$ is $\bz^d\oplus 
H-$graded.  The polynomial ring $R$ is $\bz^d$-graded by grouping all 
homogeneous elements with the same $\bz^d$-graded component.
\begin{example}
The minimal primes ideals of $I_{L}$ are $\bz^d$-graded ,  but not 
necessarily $\bz^d\oplus H-graded$.  Let $I=(x^2-y^2)\subset K[ x, y] 
$,  in this case $L=\bz(2, -2)$ and the isomorphism 
$\bz^2/L\longrightarrow \bz\oplus \bz/2\bz$,  is given by $(a, b)\mapsto 
(a+b, b \mod 2) $; it follows that deg$(x)=(1, 0)$,  deg$(y)=(1, \bar 1)$.  
The minimal primary decomposition of $I$ is given by $(x^2-y^2)=(x-y)(x+y)$ .
\end{example}
Let $A$ be any arbitrary subset of $G(S)$, we will denote by $K[ A] $ the 
$K-$vector space spanned by $A$ in $K[ G(S)] $.  If $A+S\subset A$,  we 
will call $A$ an $S-$ideal.  A proper subset $P$ of $S$ is a prime 
ideal if $P$ is an $S-$ideal and $S\setminus P$ is additively closed.  
Every $G(S)-$graded prime ideal $\underline p$ of $K[ S] $ is exactly of 
the form $K[ P] $ for some prime ideal $P$ of $S$ and the homogeneous 
localization $K[ S]_{\underline p}$ is isomorphic to $K[ S-(S\setminus 
P)] $.

\item Let $S\subset \bz^l\oplus H$ be a standard semigroup with torsion,  and 
let $I$ be a nonempty subset of $[ 1, l] $,  set $P_{I}=\{ x\in S; x_{i}> 0 { 
\rm \ for\ some\  } i\in I\}$ and $p_{I}=K[ P_{I}] $.  Then the set $\{ P_{I}\}$ is the set of 
prime ideals of $S$ (see the proof of the next lemma).  Moreover $K[ S]/p_{\{ i\}} = K[ 
S_{(i)}]$ and if $J\subset 
I$ then $P_{ J}\subset P_{ I}$.  This implies that $p_{\{ 1\}} , \ldots,  
p_{\{ l\}}$ are the unique $G(S)-$graded prime ideals of 
height one of $K[ S]$.

\end{enumerate} 
The following Lemma shows that the extension from the toric case to the lattice
 case is non trivial:
 \begin{lemma} Let ${ \cal P }\subset K[ S]$ be an $\bz^d$-graded prime 
 ideal of height $> 0$.  Then ${ \cal P }$ is $G(S)-$graded and ${ \cal P 
 }=p_{I}$ for some non empty subset $I$.
        
 \end{lemma}¥
 \demo  First,  let  remark  that if $z_{1}, z_{2}\in K[ S] $ are two pure monomials with 
  the same $\bz^d$-grade,  then $ z_{1}^h-z_{2}^h=0$,  where $h$ is the 
  order of the group $H$,  and this imply that for any prime ideal $ { 
  \bf p } $ in $K[ S] $ there exists $\xi$ a $h-$root of unity such that 
  $z_{1}+\xi z_{2}\in { \bf p }$
  
 We prove that ${ \cal P }$ contains one monomial element $t^{ 
 g_{1}}\underline u^{ g_{2}}$ for some $(g_{1}, g_{2})\in S$. 
  
  Since ht$({ \cal P })\geq 1$,  and because $I_{L}¥$ is unmixed,  ${ \cal 
  P }$ contains $Q$ an associated prime of $I_{L}¥$ and a non zero divisor $z$ for $K[ S] $. 
   Now let $z\in { \cal P }$ be a non 
  zero divisor,  we can assume that $z$ is $\bz^n-$homogeneous,  if $z$ is 
  not monomial we can write it as a sum of monomials $z=\lambda_{1}¥z_{1}+\ldots+ 
  \lambda_{r}¥z_{r}$,  with coefficients $\lambda_{i}\in K$, then for any $i=1, \ldots, r$ 
  there exists $h-$roots of unity $\xi_{i}$,  such that 
   $z_{i}+\xi_{i}¥ z_{1}\in 
  Q\subset { \cal P } $,  but $$ z=
  \sum_{{ i=1 }}^{ r } \lambda_{i}(z_{i}+\xi_{i}¥ z_{1})- (\sum_{{ i=1 
  }}^{r } \lambda_{i}\xi_{i}) z_{1}, $$ and since $z\notin Q$ we have 
  $\sum_{{ i=1 }}^{ r } \lambda_{i}\xi_{i}\not= 0$,  which implies that 
  $z_{1}\in { \cal P }$,  and we are done.
  
  The same proof shows that for  any non empty set $I$,  
    if ${ \cal P }\not\subset p_{I}$,  then we 
 can choose a monomial element $z\in { \cal P }\setminus p_{I}$.
 
 Let $I$ be the set of integers $i\in [ 1, \ldots, l]$ such that $p_{\{i\}}$ is 
  contained in ${ \cal P }$,  we will prove that $I$ is non empty and 
  $p_{I}= { \cal P }.$
  
  It is clear that $p_{I}\subset { \cal P }$,  remark that if $I=[ 
  1, \ldots, l]$,  then ${ \cal P }$ contains the unique graded maximal ideal of 
  $K[ S] $.  So we can assume that $I$ is a proper subset of $[ 
  1, \ldots, l]$.  Suppose a contrario that there exist $z \in { \cal P 
  }\setminus p_{I}$,  (If $I$ is empty choose $z$ any monomial non zero divisor),   
  we can assume that $z$ is pure monomial,  and if $z=\underline t^a 
  \underline u^b$,  then $a_{i}=0$ for all $i\in I$,  for any 
  $j\notin I$ choose a monomial $\underline t^{ c^j¥ } \underline u^{ 
  d^j }\in p_{\{ j\}}\setminus { \cal P }$,  let $c=\sum_{j\notin I} 
  c^j,  d=\sum_{j\notin I} 
  d^j$,  then $c_{j}>0$ for any $j\notin I$,  and there exist a positive 
  integer $p$ such that $p(c, d)-(a, b)\in \bn^l\oplus H \cap G=\bar S$,  
  and for some positive integer $k$ ,  $kp(c, d)-k(a, b)\in S$.  It follows 
  then that $\prod_{i\notin I}¥(\underline t^{ c^j¥ } \underline u^{ d^j 
  })^{ pk }\in{ \cal P }$ and $\underline t^{ c^j¥ } \underline u^{ 
  d^j }\in{ \cal P }$,  for some $j$. A contradiction.

\begin{theorem} Assume that the semigroup (eventually with torsion) $S $ is standard.  
Let $G(S)$  be the group generated by $S$ of rank $d$,  let 
$S':=\cap_{i=1}^l (S-(S\setminus P_i))$ be a subsemigroup of $G(S)$,  where 
$S\setminus P_i$ consist of the elements in $S$,  which the $i$ 
coordinate is 0.  Then $$K[S']=\cap_{i=1}^l K[S]_{ (p_{\{ i \}}) }, $$ 
and $K[S']$ satisfies the condition $S_{2}$. Let remark that    $K[S]_{ (p_{\{ i \}}) }$ is a 
homogeneous localization and the intersection is taken in the localization 
$T^{ -1 }K[S]$,  where $T$ is the set of all pure monomials,  also 
 since $I_{L}$ is a lattice ideal any monomial is a non zero 
divisor for $K[ S] $.  

We also have an exact sequence : $$0\la K[S] 
\la K[S'] \la K[S'\setminus S] \la 0, $$ and $ \dim K[S'\setminus S] 
\leq d-2 $.
Moreover if $S$ is simplicial then $K[S']$ is a Cohen-Macaulay ring.
\end{theorem} 
\demo It follows from \cite{GSW},  p.244,  that the 
property $S_{k}$ holds for a $\bz^d$ graded module $M$ if and only if 
$$\depth M_{(p)}\geq { \rm min \ }\{ k,  { \rm dim\ }M_{(p)}\}$$ for any 
$\bz^d$ homogeneous prime ideal $p$.

As a consequence the ring $\cap_{{\rm ht }(p)=1}¥ K[ S]_{(p)}$,  where $p$ 
runs over all $\bz^d$ homogeneous prime ideals in $ K[ S] $ of height 
one,  satisfies the condition 
$S_{2}$.  Now the above lemma proves that $\{ p_{\{ 1 \}}, \ldots, p_{\{ l 
\}}¥\}$ are all the $\bz^d$ homogeneous prime ideals in $ K[ S] $ of 
height one and then $$K[S']=\cap_{i=1}^l K[S]_{ (p_{\{ i \}}) }, $$ 
satisfies the condition $S_{2}$.  Also we have that $K[S]_{(p_{\{ i 
\}})}=K[S']_{(p_{\{ i \}})}$ and since the module $ K[S'\setminus S] $ 
is $\bz^d$-graded we get $ \dim K[S'\setminus S] \leq d-2 $.

If $S$ is simplicial,  let $x_{1}, \ldots, x_{d}$ be the variables in $R$ 
corresponding to the extreme rays of ${ \cal C }_{S}$,  then $S'$ is 
also simplicial and $x_{1}, \ldots, x_{d}$ are parameters for both $K[ 
S] ,  K[ S']$, since $S'$ satisfies the condition $S_{2}$ we have that  any pair 
$x_{i}, x_{j}$ is a regular sequence in $K[ S'] $,   if we have a relation $fx_{i}= 
\sum_{ j<i}f_{j}¥x_{j}$ then because of the grading we certainly have 
$fx_{i}= f_{j}¥x_{j}$ for some $j$,  this implies that the sequence 
$x_{1}, \ldots, x_{d}$ is a regular sequence in $K[ S'] $,  so $K[ S'] $ 
is a Cohen-Macaulay ring.

As a Corollary we have

\begin{corollary} Let $L\subset \bn^m$ be a positive lattice  of 
rank $r$,  set $\dim R/I_{L} = d=m-r$. If $\depth R/I_{L} =d-1$ then $\dim K[ S'-S] =d-2$,   $\depth D^{ d-1 }(K[ S] )=\depth K_{K[ S]} -2$,  and the following are 
equivalent:
 \begin{enumerate}
  \item $D^d(D^d(K[ S])) $ is Cohen-Macaulay.
  \item the canonical module of $K[ S] $ is Cohen-Macaulay.
  \item the module $D^{ d-1 }(K[ S] )$ is Cohen-Macaulay.
  \end{enumerate}¥
   \end{corollary}
The proof is immediate from Theorem 2.
\begin{corollary} Let $S\subset \bn^n\oplus H$ be a simplicial finitely 
generated semigroup of rank $d$,  then
 \begin{enumerate}
        \item  $D^d(D^d(K[ S])) $ is Cohen-Macaulay
  \item the canonical module of $K[ S] $ is Cohen-Macaulay
  \end{enumerate}¥

\end{corollary}

 We review the following example from \cite{K},  Example B.1: 
\begin{example} Let $K$ be a field,  $A$ the affine semigroup ring 
$$K[ 
a, b, c, d, e^2, e^3, ade, bde, cde, d^2e] .$$
We can see immediately that $K[ S'] =K[ a, b, c, d, e] $ and then it is a 
Macaulayfication of $A$,  and we get the following exact sequence (see also 
\cite{K}),  $$0\la A\la K[ a, b, c, d, e] \la C[ -1] \la 0, $$ where 
$C=A/(ad, bd, cd, d^2, e^2, e^3, ade,  bde, cde, d^2e)$ has dimension three.  
It follows that  the 
Non Cohen-Macaulay locus of $A$ is the support of $C$.
\end{example}
\begin{example}
The following example is a toric ring of codimension two and dimension 4,  which 
canonical module is not Cohen-Macaulay. The ideal $I_{L}¥\subset K[ 
a, b, c, d, e, f] =R$ has the following generators: $$ab^4c-de^3f^2, \ 
bc^3d^3-a^2e^2f^3, \  c^2d^4e-a^3b^3f, \  b^5c^4d^2-ae^5f^5, $$ $$ a^4b^7-cd^5e^4f, \  c^5d^7-a^5b^2ef^4,  
\ b^9c^5d-e^8f^7.$$

 Let $0 \la G\buildrel \phi \over \la F $ be the last term of a 
 resolution of $A:= S/I_{L}¥$,  $\sigma $ be the transpose of $\phi$,  then $ 
 F\buildrel\sigma\over \la G \la D^3(A)\la  0$ is a presentation of 
 module $D^3(A)$,  a quick computation by Macaulay gives that

$$\sigma=\pmatrix { &e^2f^2 &-bc &0 &d &-a &0 &0 &0 &0 &0\cr
 &ab^3 &-de &-f &0 &0 &c&0 &0 &0 &0\cr &c^2d^3 &-a^2f &0 &0 &0 &0 &e &-b &0 
&0 \cr &0&0 &0 &b^4c &-e^3f^2 &0 &0 &0 &-d &a \cr }$$ and that the module 
$D^3(A)$ has dimension 2,  but $\depth D^3(A)=1$.  So the canonical 
module of $S/I_{L}¥$ is not Cohen-Macaulay, in fact $\depth K_A=3$.
 
 \end{example}

 In what follows we will write $I$ instead $I_L$.
  \begin{theorem}
Let $R=K[ x_{1},  \ldots ,  x_{n}] $ be a polynomial ring,  let $A=R/I$ be a 
lattice ring,  of codimension two and dimension $d$.  If $I$ is 
minimally generated by 4 generators,  then $D^{ d-1 }(A)$ is a complete 
intersection.  In particular,  the canonical ring $K_{A}$ is 
Cohen-Macaulay of dimension d,  and the $S_{2}-$fication is a 
Macaulayfication of $A$.  The Non-Cohen-Macaulay locus of $A$ is the 
support of a Cohen-Macaulay module of dimension $d-2$.
\end{theorem}
 \demo 
  The resolution of $A$,  follows  from \cite{PS} Construction 5.2:
  
  $$0 \la R\buildrel \phi \over \la R^4 \la R^4 \la R \la A\la 0 $$ 
  where $ \sigma $ the transpose of $\phi$ is given by: 
  $$\sigma=\pmatrix{- { \underline \bf x }^{ \bf s }& { \underline \bf x 
  }^{ \bf t }&{ \underline \bf x }^{ \bf r }&-{ \underline \bf x }^{ \bf 
  p }\cr}$$ where all monomials have disjoints supports.  Then the 
  entries of $\sigma$ define a complete intersection,  that is $D^{ d-1 
  }(A)$ is a complete intersection.  The rest of the proof follows from 
  Lemma \ref{th2}.
  

{ \bf Question  }  Let $R=K[ x_{1},  \ldots ,  x_{d+2}] $ be a polynomial ring,  
let $A=R/I$ be a lattice ring  of codimension two and dimension $d$,  is 
it true that $D^{ d-1 }(A)$ has non zero divisors?

\section{ Simplicial lattices ideals of height 2}

 Let $ K $ be a field and $ 
 R:=K[y, z, x_1, \ldots, x_n]$ the ring of polynomials in the variables $ 
 y, z, x_1, \ldots, x_n$.  Let $ a_i, b_i, c_i$ $ 1\leq i\leq n $ be 
naturals 
 numbers satisfying the conditions: $$a_i\not= 0,  ( b_i, c_i )\not= 0 
 \forall i ,  ( b_1, ..., b_n)\not= 0,  ( c_1, ..., c_n)\not= 0$$

  For $i= 1, \ldots,  n$ let ${ \bf d_{i}¥ }= a_{i} { \bf e_{i}¥ }$,  
where 
  ${ \bf e_{n}¥ },  \ldots,  { \bf e_{n}¥ }$ is the canonical basis of 
  $\bn ^n$,  and  ${ \bf a_{1}¥ }=( b_1, ..., b_n),  { \bf a_{2}¥ }=( 
c_1, ..., c_n)$.  Let $H$ be a finite abelian group and 
$h_{1}, \ldots, h_{n+2}\in H $ that generates it.  Let $S$ be the 
subsemigroup of $\bn^n\oplus H$ generated by $$({ \bf d_{1}, h_{1}}) , \ldots,  
({ \bf d_{n}, h_{n}}), ({ \bf a_{1}, h_{n+1}}), ({ \bf 
a_{2}, h_{n+2}}).$$

\begin{definition} A simplicial lattice ideal  
of height two is the lattice ideal $I_{L}\subset R$,  where:
$$L=\{ { \bf w } \in \bz^{ n+2 },  w_{1}{ \bf (d_{1}, h_{1}) }+\ldots+ 
w_{n}({ \bf d_{n}, h_{n}¥ })+w_{n+1}({ \bf a_{1}, h_{n+1}})+ w_{n+2}({ \bf 
a_{2}, h_{n+2}})=0\}.$$
\end{definition}

We remark that the last two coordinates of vectors in $L$, determine 
all the lattice $L$. More precisely,   consider the group  morphism:

$$\Phi :\bz^2\longrightarrow \bz/{a_1 \bz}\times\ldots \times \bz/{a_n\bz} \ \ 
\ \
(s, p)\mapsto (sb_1-pc_1, \ldots, sb_n-pc_n)$$

The lattice $L$ is completely determined by the rank two sublattice : 
$$\tilde L \subset Ker(\Phi):=\{ (s, p)\in { \bf \bz^2}\ / \ sb_i-pc_i\equiv 0\ 
mod\ a_i,  \forall i=1, \ldots, n \}.$$
$$\tilde L =\{ (s, p)\in { \bf \bz^2}\ / \ s({ \bf a_{1}, h_{n+1}})-p({ \bf a_{2}, h_{n+2}}) \in \ 
\bz{ \bf (d_{1}, h_{1}) }+\ldots+ 
\bz({ \bf d_{n}, h_{n}¥ }) \}.$$
\begin{remark} To any vector $(s, p)\in \tilde L $ with $s\geq 0$ we 
associate a 
unique binomial $B_{(s, p)}\in I_L$ in the following way: for any
 $ i=1, \ldots, n$, let $v_i$ 
be the unique integer such that $sb_i-pc_i = v_i a_i$.  We define the 
vectors ${\bf v}_+, {\bf v}_- \in { \bn^n}$ by ${\bf v}_{+,i}= max \{ v_i, 0\} , {\bf v}_{-,i}= 
max \{- v_i, 0\}$ and we must distinguish two cases:
\begin{itemize}
\item i) if  $s\geq 0$ and  $p\geq  0$  then $B_{(s, p)}= z^s{\ul x}^{{\bf v}_{-}} - 
y^p{\ul x}^{{\bf v}_{+}}$,
\item ii)  If $s\geq  0$ and $p< 0$  then  $B_{(s, p)}= z^s y^{-p} {\ul x}^{{\bf v}_{-}} 
- {\ul x}^{{\bf v}_{+}}$.
\end{itemize}

{\bf Let $ D_i $ be the line $ D_i =\{ 
(s, p)\in \br^2 \quad \mid \quad  sb_i-pc_i =0 \}$. From now on, we suppose that 
the variables $ x_1,...,x_n $ are indexed in such a way  that the slopes of the lines $ 
D_i $ are in increasing order.}
\end{remark} 
\begin{lemma}\label{lemma6}  Consider $B= M_1-M_2\in I_{L}¥$  a 
binomial,  $ M_1, M_2$ without common factors.  We can write $B$ in only one 
of the followings forms:

\begin{enumerate}
\item $z^s-y^px_1^{ v_1}\ldots x_n^{ v_n}, \ s>0 \ v_i\geq 0\ \forall i .$

\item $y^p-z^s x_1^{ v_1}\ldots x_n^{ v_n}, \ p>0, \ v_i\geq  0 \ \forall i.$
 \item $y^pz^s-x_1^{ v_1}\ldots x_n^{ v_n}, \ p, s>0 \ v_i> 0\ \forall i.$
\item
$z^sx_1^{ v_1}\ldots x_k^{ v_k}-y^px_{ k+1}^{ v_{ k+1}}\ldots
x_{n}^{ v_n},  \ v_i\geq 0, \ p, s>0, $ and $\exists 1\leq i_1\leq k, \ 
k+1\leq i_2\leq n \ / \ v_{i_1 }, v_{i_2 }\not= 0$.  In other words if $({\bf v}, s, p) \in L$,  with $s, p>0$, and  
${ \bf v }=(v_{1}, \ldots,  v_{n})$ there exist $k $ such that 
$v_{i}<0$ 
for all $i<k $ and $v_{i}\geq 0$ for all $i\geq k $.
\end{enumerate}
\end{lemma}
As a consequence we have the following lemma:

\begin{lemma}1) There is no non trivial binomial in $I_{L}¥$,  of the type:
$ z^s y^{p} {\ul x}^{{\bf v}_{-}} - {\ul 
x}^{{\bf v}_{+}}$ with $s \geq 0$,  $p\geq 0$,  and ${{\bf v}_{-}} \not=0$.

\noindent 2) Consider an equality (where every fraction is reduced): $${{z^{s} P_{1}¥(\underline x)}\over{ 
\underline x^{ \bf v}}}={{y^{p} P_{2}¥(\underline x)}\over{ \underline 
x^{\bf w}}}.$$ 
\begin{itemize}
\item If there exist an index 
$i_{1}$ such that $v_{ i_{1} }>0,  w_{ i_{1} }=0 $,  then for all 
$j\geq i_{1}$,  $v_{j}> w_{j} $. 
\item If there exist an index 
$i_{2}$ such that $v_{ i_{2} }=0,  w_{ i_{2} }>0 $,  then for all 
$j\leq i_{2}$,  $v_{j}< w_{j} $.
\end{itemize}
\end{lemma}
 The following proposition is an extension of \cite{Mo1}, to the lattice case. 
 \begin{proposition}  We can describe a fan decomposition 
of 
${ \bf R_+^2}$,  more precisely we have vectors $\varepsilon_{ -1},  
\varepsilon_0, ..., \varepsilon_{ m+1}\in \tilde L \cap { \bf Z_+^2}$ 
such 
that 
\begin{itemize}
\item $\varepsilon_{ -1}=(s_{-1} , 0) ,  
\varepsilon_0 =(s_0,  p_{0})$,  with $0\leq s_{0}\leq s_{-1}$.
\item Consider the 
Euclidean algorithm to compute the $gcd(s_{-1}, s_{ 0})$ :
 \begin{eqnarray*}
        s_{-1} &=& q_1s_0-s_1 \cr s_0 &= &
q_2s_1-s_2 \cr & \ldots &\cr s_{m-1} &=& q_{m+1}s_m \cr s_{m+1} &=&0 
\cr 
\end{eqnarray*} 
        $$q_i\ge 2 \ , \ s_i\ge 0\ \ \forall i $$
Let $p_i$ be the sequence of integers defined by
  \begin{eqnarray*}
  p_{ i+2}&=&q_{ i+2}p_{i+1}-p_{ i}\ \ \ \ -1\leq i\leq 
 m-1
 \end{eqnarray*} then $\varepsilon_i =(s_i,  p_{i})$.
 \item ¥$ \varepsilon_i, \varepsilon_{ i+1} $ is a basis of $ \tilde L  $ 
and $det(\varepsilon_i, \varepsilon_{ i+1})= p_{0}s_{-1}>0$.  
\end{itemize}¥
\end{proposition}

Note that the existence of the basis $\varepsilon_{ -1}, \varepsilon_0$ is provided by \cite{Co}, page 62.
\begin{definition}
 Let $r_{ j, i}$ be the sequence of integers defined by
  \begin{eqnarray*}
   r_{ j, i}&=&(s_ib_j-p_ic_j)/a_j \ \ \ \ -1\leq i\leq m+1 , \ 1\leq 
 j\leq n
  \end{eqnarray*}
and $ { \bf r}_{i} $ the vector with coordinates $ r_{ j, i}$.
\end{definition}

\begin{lemma}
1) Any of the sequences $s_i,  p_i, r_{ j, i},  \ 1\leq j\leq n $ satisfy 
the recurrent relation:

 $$v_{ i+2}=q_{ i+2}v_{i+1}-v_{ i}{ \rm \  for \ }-1\leq i\leq m-1.$$

2) The sequences  $ s_i, r_{ j, i}$ (for all $  j $)
are  strictly  decreasing but the sequence $ p_i$ 
is  strictly  increasing. 

3) Set   $ \nu $  (resp. $ \mu$ )
the greatest integer $j$ such that $ { \bf r }_{j }={ \bf r }_{j, + }$ 
(resp. the 
smallest integer $j$ such that $ { \bf r }_{j }=- { \bf r }_{j, - }$),  
then $-1\leq 
\nu\leq \mu\leq m .$

4) $\supp { \bf r}_{i+1, +}\subset \supp { \bf r}_{i, +}$.
\end{lemma}

\begin{theorem}

1) The ring $R/I_L$ is  arithmetically Cohen-Macaulay 
if and only if $ \mu=\nu $. In
this case the ideal 
 $ I_L$  is generated by: 

\begin{eqnarray*}
F&= & z^{s_\nu }\ -\ y^{p_\nu } {\underline x}^{{ \bf r}_{\nu}}\cr 
G&=& 
y^{p_{\nu +1}} \ -\ z^{s_{\nu +1}}{\underline x}^{{ -\bf r}_{\nu+1}} 
\cr 
H&=& z^{s_{\nu }-s_{\nu +1}}y^{p_{\nu+1 }-p_{\nu }}
{\underline x}^{{ \bf r}_{\nu}\ - { \bf r}_{\nu+1}}\cr
\end{eqnarray*}¥

2) If $R/I_L$   is not  arithmetically Cohen-Macaulay
the ideal  $I_L$
is
 generated by  $\tau := 3+(q_{\nu +2}-1)+...+(q_{\mu+1}-1)$
 equations:
 \begin{eqnarray*}
 z^{s_\nu } &-&y^{p_\nu } {\underline x}^{{ \bf r}_{\nu} } \cr 
 y^{p_{\nu +1}-p_{\nu } } z^{s_{\nu }-s_{\nu +1}}&-& {\underline 
 x}^{{ \bf r}_{\nu }-{ \bf r}_{\nu +1}}\cr z^{s_{\nu +1}}{\underline 
x}^{{ \bf r}_{\nu 
 +1, -}} &-&y^{p_{\nu +1}}{\underline x}^{{ \bf r}_{\nu +1, +}} \cr 
 y^{2p_{\nu+1}- p_{\nu }} z^{s_{\nu }-2s_{\nu +1}}&-& {\underline 
 x}^{{ \bf r}_{\nu }-2{ \bf r}_{\nu+1 }}\cr &\ldots& \cr y^{(q_{\nu 
+2}-1)p_{\nu 
 +1}-p_{\nu }} z^{s_{\nu }-(q_{\nu +2}-1)s_{\nu +1}} &-& {\underline 
 x}^{{ \bf r}_{\nu }-(q_{\nu +2}-1){ \bf r}_{\nu +1}}\cr z^{s_{\nu 
+2}}{\underline 
 x}^{{ \bf r}_{\nu +2, -}} &-&y^{p_{\nu +2}}{\underline x}^{{ \bf 
r}_{\nu +2, +}}\ \cr 
 &\ldots& \cr &\ldots& \cr z^{s_{\mu }}{\underline x}^{{ \bf r}_{\mu 
, -}} 
 &-&y^{p_{\mu }}{\underline x}^{{ \bf r}_{\mu, + }} \cr 
y^{2p_{\mu}-p_{\mu-1 }} 
 z^{s_{\mu-1}-2s_{\mu}}&-& {\underline x}^{{ \bf r}_{\mu-1}-2{ \bf 
r}_{\mu}}\cr 
 &\ldots& \cr y^{(q_{\mu+1}-1)p_{\mu}-p_{\mu-1 }} 
 z^{s_{\mu-1}-(q_{\mu+1}-1)s_{\mu}} &-&{\underline x}^{{ \bf 
r}_{\mu-1 
 }-(q_{\mu+1}-1){ \bf r}_{\mu}}\cr z^{s_{\mu+1}}{\underline 
 x}^{{ \bf r}_{\mu+1, -}}&-&y^{p_{\mu+1}}\cr
 \end{eqnarray*}¥

They form a Groebner's basis for the reverse lexicographic order with 
respect to $z<y<x_1<\ldots <x_n $.
\end{theorem}
 \demo Note that the proof given in \cite{Mo1}, pp.1089, applies here without restriction. 
We  outline the proof of 2): it consist to prove that the leading term of any binomial in $I_L$
 for the reverse lexicographic order with 
respect to $z<y<x_1<\ldots <x_n $ is a factor of the leading term of some binomial in the above list. For example, let
$B$ be a binomial corresponding to the lattice point $({\bf v}, s,p)$ with $p\geq 0, s\geq 0$. By the fan decomposition of
 $\br_+\times\br_+$, there exists some $i\geq -1$ such that $(p,s)= \alpha \varepsilon_{ i}+  
\beta \varepsilon_{ i+1},$ with intehers $\alpha > 0 , \beta \geq 0$, this  imply ${\bf v}=
\alpha {\bf r}_{ i}+  \beta {\bf r}_{ i+1}$. We need to consider three cases:
\begin{itemize}
\item if $i<\nu$ then the coordinates of  ${\bf r}_{ i},{\bf r}_{ i+1}$ are all positive, so the leading term of $B$ is $z^s$ 
but $s=\alpha {s}_{ i}+  \beta {s}_{ i+1}\geq {s}_{ \nu }.$
\item By a similar argument if $i\geq \mu+1$ then the coordinates of  ${\bf r}_{ i},{\bf r}_{ i+1}$ are all negative,
 so the leading term of $B$ is $y^p$ 
but $p=\alpha {p}_{ i}+  \beta {p}_{ i+1}\geq {p}_{ \mu+1 }.$
\item if $\nu \leq i\leq \mu $ then  the leading term of $B$ is $z^s{\underline x}^{{ \bf v}_{-}} $ which is a factor of 
$z^{s_{i }}{\underline x}^{{ \bf r}_{i 
, -}} $.
\end{itemize}
If 
$B$ is a binomial corresponding to the lattice point $({\bf v}, s,p)$ with $p< 0, s\geq 0$.
 We argue with similar arguments
using  the fan decomposition of
 $\br_+\times\br_-$ (that is every two consecutive vectors is a basis of $\tilde L$),
 given by the sequence of vectors $$\varepsilon_{ -1}-\varepsilon_{ 0},...,
\varepsilon_{ -1}-(q_1-1)\varepsilon_{ 0}=\varepsilon_{ 0}-\varepsilon_{ 1},...,
\varepsilon_{ 0}-(q_2-1)\varepsilon_{ 1}=\varepsilon_{ 1}-\varepsilon_{ 2},...,$$
$$
\varepsilon_{ m-1}-(q_m-1)\varepsilon_{ m}=\varepsilon_{ m}-\varepsilon_{ m+1},-\varepsilon_{ m+1}.$$

\section{ Macaulayfication of codimension two simplicial toric rings}
The aim of this section consist to give an explicit description of the 
semigroup $S'$ such that $K[S']$ is the Macaulayfication of the 
simplicial semigroup ring of codimension two $K[S]$.  (see Theorem 5):
 
 We recall that  $S'=\cap_{i=1}^l (S-(S\setminus P_i))$ is a subsemigroup of 
 $G(S)$,  where $S\setminus P_i$ consist of the elements in $S$,  which 
 the $i$ coordinate is 0.  the ring $$K[S']=\cap_{i=1}^l K[S]_{ (p_{\{ 
 i \}}) }, $$ is a Cohen-Macaulay ring,  where $K[S]_{ (p_{\{ i \}}) }$ 
 is a homogeneous localization and the intersection is in the 
 localization $T^{ -1 }K[S]$,  where $T$ is the set of all pure 
 monomials.  Remark that since $I_{L}$ is a lattice ideal any monomial 
 is a non zero divisor for $K[ S] $.
Any simplicial group is trivially standard. In what follows we will write
 $I$ instead of $I_L$.

\begin{remark} Since $S$ is simplicial of codimension two,  every 
element in $S-(S\setminus P_i)$ can be  viewed as a quotient of monomials 
$\displaystyle {{M(y, z, \underline x)}\over {N(y, z, \underline x)}}$ where 
$M, N$ are monomials with disjoints supports and $N\notin p_{\{ i \}}$,  we 
notice that

\begin{equation}
p_{\{ i \}}=\cases { (x_{i}, y, z)& if $b_{i}\not= 0$ and $c_{i}\not= 0$,\cr
 (x_{i},  y)& if $b_{i}= 0$ and $c_{i}\not= 0$,\cr
 (x_{i}, z)& if $b_{i}\not= 0$ and $c_{i}= 0$.\cr}
\end{equation}¥
 
\end{remark} 
\begin{lemma} Let $E\in \cap_{i=1}^n (S-(S\setminus P_i))$, then for 
each 
$i$ we can write $\displaystyle E= 
{{z^{\alpha_i}y^{\beta_i}P_i(\underline x)}\over{Q_i(\underline x 
)}}$,  such that $x_i$ is not in the support of $Q_i$.
\end{lemma}

\demo The assertion is clear if $b_{i}\not= 0$ and $c_{i}\not= 0$ for all $i=1, ..., n$.
 If $b_{i}= 0$ and $c_{i}\not= 0$ then $s_{j}b_i-p_{j}c_i=-p_{j}c_i\leq 0$ and we have equality
 if and only if $i=-1$,  this implies $\nu=-1,  p_{\nu}=0$. 
 Regarding the order introduced in  the variables  $x_1, ..., x_n$ by the lemma 3,  
we can suppose that there exist natural integers $k, l$ such that 
$b_{1}=\ldots, b_{k}=0,  a_{n-l}=\ldots a_{n}=0$ and that $k, l$ are the 
biggest possible.
It will be enough to prove the Lemma for  $i\leq k $ and  $0<k<n$. 
 Let  $\displaystyle E\in \cap_{i=1}^n (S-(S\setminus P_i)),  
 E= {{y^{\beta_i}P_i(\underline 
x)}\over{Q_i(\underline x ){z^{\alpha_i}}}}$ where  $x_i$ is not in the 
support of $Q_i$, and $\alpha_i>0$. On the other hand for  any $k< j<n-l$ we have 
$\displaystyle E= {{z^{\alpha_j}y^{\beta_j}P_j(\underline 
x)}\over{Q_j(\underline x )}}$,  where $x_j$ is not in the support of 
$Q_j$ remark that $x_i$ belongs to the support of $Q_j$, otherwise we have finish our proof,  we can also  assume that $P_j $ and 
$Q_j$ have disjoint support,  this gives us the following element in 
$I$:

$$y^{\beta_i}P_i(\underline 
x)Q_j(\underline x )-z^{\alpha_j+\alpha_i}y^{\beta_j}P_j(\underline 
x)Q_i(\underline x )$$

Since $x_{i}$ appears in the left side of this equality but no in the 
right side,  we must have $\beta_j>\beta_i$.  More precisely we write 
$Q_{j}(\underline x )=x_{i}^{\gamma_{i}}\tilde Q_{j}(\underline x )$,  
$P_{i}(\underline x )=x_{i}^{\delta_{i}}\tilde P_{i}(\underline x )$,  
and since the couple $(\alpha_i+\alpha_j, \beta_j- \beta_i)$ belongs to 
the lattice $\tilde L$,  there exist integers $A, B$ such that:
$$(\alpha_i+\alpha_j, \beta_j- \beta_i)=A(s_{-1}, 0 )+B(s_{0}, p _{0} )$$
 this implies that $\beta_j- \beta_i= Bp _{0} $. We have the following elements in $I$
 $$z^{s_{0}}{\underline 
x}^{{ \bf r}_{0, -}} -y^{p_{0}}{\underline x}^{{ \bf r}_{0, +}}, 
 z^{Bs_{0}}{\underline 
x}^{{ B\bf r}_{0, -}} -y^{Bp_{0}}{\underline x}^{{B \bf r}_{0, +}} $$
this implies $\gamma_i+\delta_i=B{ \bf r}_{0, -, i}$, 
 and we have the following equality:
 $${{z^{Bs_{0}}{\underline x}^{{B \bf \hat r}_{0, -}}{x^{\delta_i}}}\over
{{\underline x}^{{B \bf r}_{0, +}}}}={{y^{Bp_{0}}}\over {{ x_i}^{\gamma_i}}}=
{{{y^{\beta_j- \beta_i}}}\over {{ x_i}^{\gamma_i}}}$$
where we have set ${ \bf r}_{0, -, i}$ for the $i-$coordinate of the vector
 ${ \bf r}_{0, -}$ and ${ \bf \hat r}_{0, -}$is the vector ${ \bf r}_{0, -}$   with the $i-$coordinate equal to zero.
Finally we have $$E= {{z^{\alpha_j}y^{\beta_i}P_i(\underline 
x)}\over{\tilde Q_j(\underline x )}}\times {{z^{Bs_{0}}{\underline x}^{{B \bf \hat r}_{0, -}}{x^{\delta_i}}}\over
{{\underline x}^{{B \bf r}_{0, +}}}}$$
and $x_i$ is not in the support of the denominator,  and we are done. 
 \begin{lemma} Let 
$E\in \cap_{i=1}^n (S-(S\setminus P_i))$.  For each $i$ we write 
$\displaystyle E= {{z^{\alpha_i}y^{\beta_i}P_i(\underline 
x)}\over{Q_i(\underline x )}}$,  where $x_i$ is not in the support of 
$Q_i$,  and we can assume that $P_i $ and $Q_i$ have disjoint support.  
Then:
\begin{enumerate}
\item For all $i$,  we can assume that $\alpha_i<s_\nu$ and 
$\beta_i<p_{\mu+1}$. The equality 
$\displaystyle{{y^{\beta_i}P_i(\underline x)}\over{Q_i(\underline x 
)}}={{y^{\beta_j}P_j(\underline x)}\over{Q_j(\underline x )}}$,  
such that $x_i$ is not in the support of $Q_i$ and $x_j$
 is not in the support of $Q_j$,  implies $i=j$ and this 
equality is an identity.  The same is true for $z$.
\item If there exists some index $i$ such that  $\alpha_i=\beta_i=0$ 
then $E\in S$.
\item We can write  $E=z^{\alpha}y^{\beta}E' $ where $E'\in S'$,  
where $\alpha$ is the minimum of all the $\alpha_i$ and $\beta$ 
is the minimum of all the $\beta_i$. In particular we can assume that 
there exist indexes $i, j$ such that $\alpha_i=0$ and $\beta_j=0$ .

\item If $\displaystyle E= 
{{z^{\alpha}P(\underline x)}\over{Q(\underline x )}}$,  where $P $ and 
$Q$ have disjoint support,  then $x_1$ is not in the support of $Q.$

\end{enumerate}

 \end{lemma}
 \demo
 \begin{enumerate}
\item Suppose that  $\displaystyle{{y^{\beta_i}P_i(\underline 
x)}\over{Q_i(\underline x )}}={{y^{\beta_j}P_j(\underline 
x)}\over{Q_j(\underline x )}}$ such that $x_i$ is not in the support 
of $Q_i$,  $x_j$
 is not in the support of $Q_j$,  $i\not= j$ and $\beta_j \geq \beta_i $.  It 
 follows that ${ \tilde P}_i(\underline x){ \tilde 
 Q}_j(\underline x )-y^{\beta_j-\beta_i}{ \tilde P}_j(\underline x){ \tilde 
 Q}_i(\underline x )$ belongs to $I$ where ${ \tilde P}_{l}¥=P_{l}/\gcd(P_{i}, P_{j})$,  ${ 
 \tilde Q}_{l}¥=Q_{l}/\gcd(Q_{i}, Q_{j})$.  Since $0\leq \beta_j-\beta_i < 
 p_{\mu+1}¥$ such element cannot 
 exists in $I$,  and we are done. 
  \item Suppose that $\displaystyle E= 
 {{P_i(\underline x)}\over{Q_i(\underline x )}}$,  such that $x_i$ is 
 not in the support of $Q_i$ but $Q_i\not= 1$ and $P_i, Q_i$ have 
 disjoint support.  Let $x_{j}$ be in the support of $Q_i$,  then we can 
 write $\displaystyle E= {{z^{\alpha_j}y^{\beta_j}P_j(\underline 
 x)}\over{Q_j(\underline x )}}$,  such that $x_j$ is not in the support 
 of $Q_j$.  It follows that $z^{\alpha_j}y^{\beta_j}P_j(\underline 
 x)Q_i(\underline x )- P_i(\underline x)Q_j(\underline x )$ belongs to 
 $I$. We get  a contradiction since $x_j$ is not 
 in the support of $P_iQ_j$. 
  \item It is clear that 
 $\displaystyle {{z^{\alpha_i-\alpha }y^{\beta_i-\beta }P_i(\underline 
 x)}\over{Q_i(\underline x )}}= {{z^{\alpha_j-\alpha }y^{\beta_j-\beta 
 }P_j(\underline x)}\over{Q_j(\underline x )}}$ in the field of 
 fractions of $K[ S] $.  We set $E'=\displaystyle {{z^{\alpha_i-\alpha 
 }y^{\beta_i-\beta }P_i(\underline x)}\over{Q_i(\underline x )}}$,  now 
 it is clear that $E'\in S'$ and $E=z^\alpha y^\beta E'$.  
 \item Suppose that $x_1$ is in in the support of $Q$,  then 
we can write $ \displaystyle  {{z^{\alpha}P(\underline x)}\over{Q(\underline x 
)}}={{z^{\alpha_1}y^{\beta_1}P_1(\underline x)}\over{Q_1(\underline x 
)}}$,  such that $x_1$ is not in the support of $Q_1$ and $\beta_1>0$.  It follows then that 

$z^{\alpha}P(\underline x) Q_1(\underline x )- z^{\alpha_1}y^{\beta_1}P_1(\underline x)
Q(\underline x ) \in I$,  but lemma 4,  implies that $\alpha > 
\alpha_1$,  and we get that $z^{\alpha-\alpha_1}P(\underline x) Q_1(\underline x )- y^{\beta_1}P_1(\underline x)
Q(\underline x ) \in I$,  and $x_{1}$ is the support of $P_1(\underline x)
Q(\underline x )$ but not in the support of $P(\underline x) Q_1(\underline 
x )$,  applying again lemma 4,  we get a contradiction since 
$\alpha-\alpha_1<s_{\nu}$.

 \end{enumerate}
 
\begin{theorem} 
\begin{enumerate}

\item Any element in the minimal basis in $I$ of the type
 $$z^{s_{\nu+l}}\underline x^{{ \bf r }_{\nu+l}, -} 
-y^{p_{\nu+l}}\underline 
 x^{{ \bf r }_{\nu+l}, +}, $$ for $1\leq l\leq \mu-\nu$,  gives rise to 
a non 
 trivial element $$\displaystyle E_l={{y^{p_{\nu+l}}}\over{\underline 
 x^{{ \bf r }_{\nu+l}, -}}}= {{{z^{s_{\nu+l}}}\over{\underline 
 x^{{ \bf r }_{\nu+l}, +}}}}\in S' .$$ 
   \item Any element $E\in S'$ 
 which can be written as $$\displaystyle{{y^{\beta}}\over{\underline 
 x^{{\bf v}_-}}}= {{z^{\alpha}}\over{\underline x^{{\bf v}_+}}}, $$ where 
${\bf v}_+, {\bf v}_-$ 
 have disjoint support,  belongs to the semigroup generated by $S$ and 
 the elements $E_l$,  for $1\leq l\leq \mu-\nu$.
  \item Any element $E\in S'$ which can be written 
 as $$\displaystyle {{y^{\beta} P(\underline x))}\over{\underline 
 x^{{\bf v}_-}}}= {{z^{\alpha}Q(\underline x)}\over{\underline x^{{\bf v}_+}}}, $$ 
 where ${\bf v}_+, {\bf v}_-$ have disjoint support, belongs to the semigroup 
 generated by $S$ and the elements $E_l$ for $1\leq l\leq \mu-\nu$.

\end{enumerate}
 \end{theorem}
 \demo
 \begin{enumerate}

\item It is clear that $\displaystyle 
E_l={{y^{p_{\nu+l}}}\over{\underline x^{{ \bf r }_{\nu+l}, -}}}  
{{{z^{s_{\nu+l}}}\over{\underline x^{{ \bf r }_{\nu+l}, +}}}}\in S'$.  
We have 
  $E_l\notin S$ since $p_{\nu+l}<p_{\mu+1}$. 
   \item Let $E\in S'$ such that 
  $\displaystyle{{y^{\beta}}\over{\underline x^{{\bf v}_-}}}= 
  {{z^{\alpha}}\over{\underline x^{{\bf v}_+}}}$,  where ${\bf v}_+, {\bf v}_-$ have 
  disjoint support.  It follows that $y^{\beta}{\underline x^{{\bf v}_+}}- 
  z^{\alpha}{\underline x^{{\bf v}_-}}$ belongs to $I$ then $(\alpha , \beta 
  )\in \ker \Phi $ and there exist positive integers $k,   
  \lambda_{1}, \lambda_{2}$ such that $$(\alpha , \beta )= 
  \lambda_{1}(s_{k},  p_{k})+ \lambda_{2}(s_{k+1},  p_{k+1})$$ and as 
  consequence of this $$v_{j}=\lambda_{1}r_{j, k}+ 
\lambda_{2}r_{j, k+1} \ 
  { \rm \ for \ all\ }1\leq j\leq n. $$  We recall that if 
$r_{j, l}>0$ 
  for some $j, l$ then $r_{m, l}>0$ for all $m>j$ and that $\supp 
  { \bf r }_{k+1, +}\subset \supp { \bf r }_{k, +}$.  By lemma 3, 
  there exist $\delta $ such that $v_{l}<0$ for all $l<\delta $ and 
  $v_{l}\geq 0$ for all $l\geq \delta $.  Let $\supp { \bf r 
}_{k, +}=\{ 
  l_{1}, \ldots, n \},\  \supp { \bf r }_{k+1, +}=\{ l_{2}, \ldots, n \}$,  
with 
  $l_{1}\leq l_{2}$.  It follows that $ l_{1} \leq \delta \leq l_{2}$ 
  then we can write $$E= {{z^{\alpha}}\over{\underline x^{{\bf v}_+}}}= 
  ({{{x_{l_{1}¥ }^{r_{l_{1} , k}}\ldots x_{\delta -1¥}^{r_{\delta 
  -1, k}}z^{s_{k}}}\over{\underline x^{r_{k}, +}}}})^{ \lambda_{1}¥ } 
  ({{{x_{\delta }^{-r_{\delta , k+1}}\ldots 
  x_{l_{2}-1¥}^{-r_{l_{2}-1, k+1}}¥z^{s_{k+1}}}\over{\underline 
  x^{r_{k+1}, +}}}})^{ \lambda_{2}¥ }$$ so $$E=(x_{l_{1}¥ }^{r_{l_{1} 
  , k}}\ldots x_{\delta -1¥}^{r_{\delta -1, k}}x_{\delta }^{-r_{\delta 
  , k+1}}\ldots x_{l_{2}-1¥}^{-r_{l_{2}-1, k+1}})^{ \lambda_{2}¥ } 
  ¥E_{k}^{ \lambda_{1}¥ } E_{k+1}^{ \lambda_{2}¥ }.$$
  
  \item If $\displaystyle 
 E={{y^{\beta} P(\underline x)}\over{\underline x^{{\bf v}_-}}}= 
 {{z^{\alpha}Q(\underline x)}\over{\underline x^{{\bf v}_+}}}$,  where 
 ${\bf v}_+, {\bf v}_-$ have disjoint support,  then after division by the common 
 factor of $P,  Q$ we can assume that they have disjoint support.  But then 
 $P(\underline x){\underline x^{{\bf v}_+}}$ and $Q(\underline 
x){\underline x^{{\bf v}_-}}$ 
 have disjoint support.
 It follows that the element $\displaystyle 
 E':={{y^{\beta} }\over{Q(\underline x){\underline x^{{\bf v}_-}}}}= 
 {{z^{\alpha}}\over{P(\underline x){\underline x^{{\bf v}_+}}}}$ belongs to $S'$
 and we can write $E= P(\underline x)Q(\underline x)E'$, the assertion follows 
 from the previous item.
 \end{enumerate}
 \begin{theorem}

 Any element $E\in S'$ belongs to 
 the semigroup generated by $S$ and the elements $E_l$ for $1\leq l\leq \mu-\nu$.

 \end{theorem}

\demo Let $E\in S' $ be a non trivial element,  by lemma 7, item 4 we can 
write 
$\displaystyle 
 E={{z^{\alpha} P_{1}¥(\underline x)}\over{ x_{1}^{a_1^1}\ldots 
 x_{n}^{a_n^1}}}$ with $a_1^1=0 $.  Let $i_{1}$ be the biggest integer  such 
 that $a_j^1=0 $ for $j< i_{1}$ but $a_{ i_{1} }^1>0$.  Since $E\in S'=\cap_{i=1}^n (S-(S\setminus P_i))$  we 
 can write $E={{y^{\beta_{1}¥}z^{\alpha_{1}¥} P_{3}¥(\underline 
 x)}\over{ x_{1}^{a_1^3}\ldots x_{n}^{a_n^3}}}$ with $a_{ i_{1} }^3=0$ 
 and $\beta_{1}>0$.  It then follows that 
 $${z^{\alpha} P_{1}¥(\underline x)x_{1}^{a_1^3}\ldots x_{n}^{a_n^3}-
 {y^{\beta_{1}¥}z^{\alpha_{1}¥} P_{3}¥(\underline 
 x)}{ x_{1}^{a_1^1}\ldots x_{n}^{a_n^1}}}\in I,$$ and lemma 4 implies that  
 $0< \alpha_{1}<\alpha$,  so we have that 
 $${z^{\alpha-\alpha_{1}} P_{1}¥(\underline x)x_{1}^{a_1^3}\ldots x_{n}^{a_n^3}-
 {y^{\beta_{1}¥} P_{3}¥(\underline 
 x)}{ x_{1}^{a_1^1}\ldots x_{n}^{a_n^1}}}\in I.$$ 
 Since $a_{ i_{1} }^1>0,  a_{ i_{1} }^3=0 $ then for all 
 $j\geq i_{1}$,  $a_{ j}^1\geq a_{j}^3 $ by lemma 4. 
 
Thus we can write the equality:
\begin{equation}{{z^{\alpha -\alpha_{1}} P_{1}¥(\underline x)}\over{ 
x_{1}^{a_{1}^1}\ldots 
x_{i_{1}}^{a_{i_{1}}^1}x_{i_{1}+1}^{a_{i_{1}+1}^1-a_{i_{1}+1}^3}\ldots 
x_{n}^{a_n^1-a_n^3}}}={{y^{\beta_{1}¥} P_{3}¥(\underline x)}\over{ 
x_{1}^{a_1^3}\ldots x_{k_{1}-1}^{a_{k_{1}-1}^3}}} 
\end{equation} 
Since the denominators have disjoint support,  this equality gives one element $F_{1}\in S'$ that  belongs 
to the semigroup generated by $S$ and $E_{1}, \ldots,  E_{\mu-\nu}$.  
Then we have that:
\begin{equation} 
E=F_{1}{{z^{\alpha_{1}}} \over{ x_{i_{1}+1}^{a_{i_{1}+1}^3}\ldots 
x_{n}^{a_n^3}}} 
\end{equation}
Now either $a_j^3=0$ for all $j>i_{1}$,  and in this case we have 
finished the proof of the theorem,  or there exist $i_{2}> i_{1}$ 
such that $a_j^3=0$ for all $ i_{1}\leq j<i_{2}$,  but $a_{i_{2}}^3>0$.
 Since $E\in S'=\cap_{i=1}^n (S-(S\setminus P_i))$  we 
 can write $E=\displaystyle {{y^{\beta_{2}¥}z^{\alpha_{2}¥} P_{4}¥(\underline 
 x)}\over{ x_{1}^{a_1^4}\ldots x_{n}^{a_n^4}}}$ with $a_{ i_{2} }^4=0$ 
 and $\beta_{2}\geq 0$. We have the following element in $I_L$
   
 $$y^{\beta_2}z^{\alpha_2} P_{4}¥(\underline x){\underline x}^
{({{\bf a}_3}-{\bf a}_4)_+}- y^{\beta_1}z^{\alpha_1} P_{3}¥(\underline x){\underline x}^
{({{\bf a}_3}-{\bf a}_4)_-}$$
 First since $\alpha_{1}< s_\nu, \alpha_{2}< s_\nu,  \beta_{1}<p_{\mu +1}, \beta_{2}<p_{\mu +1} $,  we must have $\alpha_{1}\not=\alpha_{2},  \beta_{1}\not=\beta_{2}. $

 Now suppose that  $\beta_{2}< \beta_{1}$,  we have two cases:
 \begin{enumerate}
 \item If $\alpha_{1}>\alpha_{2}$ then Lemma 4,  implies that $a_i^3>a_i^4$ for all $i$ but $a_{i_1}^3=0$, 
 this is a contradiction.
\item If $\alpha_{1}<\alpha_{2}$ since $a_{i_2}^3>0, a_{i_2}^4=0$ by Lemma 4, we get  $a_i^3>a_i^4$ 
for all $i\leq i_2$ but $a_{i_1}^3=0$,  this is a contradiction. 

 \end{enumerate}
 So we have  $\beta_{2}>\beta_{1}$,  if we assume that $\alpha_{1}<\alpha_{2}$ since $a_{i_2}^3>0, a_{i_2}^4=0$ 
by lemma 4 we have  $a_i^3>a_i^4$ for all $i\leq i_2$ but $a_{i_1}^3=0$,  this is a contradiction. 
Finally we get  $\beta_{2}>\beta_{1}$ and $\alpha_{1}>\alpha_{2}$.

 Using lemma 4,  we argue as before and we get that for all 
 $j\geq i_{2}$,  $a_{ j}^3\geq a_{j}^4 $.  We have the following equality:
$${{z^{\alpha_{1}-\alpha_{2}}} \over{ x_{i_{2}}^{a_{i_{2}}^3-a_{i_{2}}^4}\ldots 
x_{n}^{a_n^3-a_n^4}}}={{y^{\beta_{2}¥ -\beta_{1}} P_{4}¥(\underline x)} 
\over{ 
 P_{{ 3}}¥(\underline x) x_{1}^{a_1^4} \ldots 
 x_{i_{2}-1}^{a_{i_{2}-1}^4}}}.$$
This equality defines one element $F_{2}\in S'$ that  belongs 
to the semigroup generated by $S$ and $E_{1}, \ldots,  E_{\mu-\nu}$,  and 
we have

\begin{equation} 
E=F_{1}F_{2}{{z^{\alpha_{2}}} \over{ x_{i_{2}+1}^{a_{i_{2}+1}^{ 
4 }}\ldots x_{n}^{a_n^{ 4 }}}}.
 \end{equation}
We can continue and we can write $$E=F_{1}F_{2}\ldots F_{m},$$ where 
$F_{1}, F_{2}, \ldots, F_{m}$ belong to the semigroup generated by 
$S$ and $E_{1}, \ldots,  E_{\mu-\nu}$. This ends the proof of the theorem.

 \begin{example}
 Let $k$ be a non zero natural number,  and consider the simplicial toric 
 variety defined parametrically by:
 $$x_{1}=u_{1}^{ 2k }, \ldots,  x_{k}=u_{k}^{ 2k },  y= u_{1}^{ k+1 
 }u_{2}u_{3}¥\ldots u_{k},  z= u_{1}u_{2}^{ k+1  }u_{3}¥\ldots u_{k}$$
 It is a codimension two variety in $\bp^{ k+1 }$.  Let $I_{k} $ be the 
 vanishing ideal of this variety.  We apply the algorithm described in 
 proposition 1 to find a system of generators of $I_{k} $:
 
 \begin{eqnarray*}
         y^{ 2k }& - & x_{1}^{ k+1 } x_{2} x_{3}\ldots x_{k} \\
 z^2x_{1} & - & y^2x_{2} \\
 y^{ 2k-2 } z^2& - & x_{1}^{ k } x_{2}^2x_{3}\ldots x_{k} \\
 y^{ 2k-4 } z^4& - & x_{1}^{ k-1 } x_{2}^3x_{3}\ldots x_{k} \\
 ¥ & \ldots & ¥\\
 y^2 z^{ 2k-2}& - & x_{1}^2 x_{2}^{ k }x_{3}\ldots x_{k} \\
 z^{ 2k}& - & x_{1} x_{2}^{ k +1}x_{3}\ldots x_{k}
 \end{eqnarray*}¥
 In order to get the Macaulayfication we must consider the element:
 $${{y^2}\over {x_{1}}}={{z^2}\over {x_{2}}}= u_{1}^2 u_{2}^2u_{3}^2\ldots u_{k}^2.$$
 The Macaulayfication will be the semigroup ring :
 $$K[ S'] = K[u_{1}^{ 2k }, \ldots, u_{k}^{ 2k },   u_{1}^{ k+1 
 }u_{2}u_{3}¥\ldots u_{k},  u_{1}u_{2}^{ k+1  }
 u_{3}¥\ldots u_{k}, u_{1}^2 u_{2}^2u_{3}^2\ldots u_{k}^2]  .$$
 In fact it is easy to check that 
 $$ K[ S']= K[ x_{1},  x_{2},  x_{3},  \ldots,  x_{k}, y,  z,  w] / (z^2-x_{2}w,  
 y^2-x_{1}w,  w^k- x_{1} x_{2} x_{3} \ldots x_ {k}), $$ and it is a 
 complete intersection.
 \end{example}
 \begin{example} We  can apply our methods to some non toric cases. 
   The (non-toric) variety $V\subset \Bbb A^7$   defined by
$$x_1=s^4+t^4;\ x_2=s^2tu;\ x_3=s^3t;\ x_4=st^3;\ x_5=su^3;\
x_6=s^2t^2v;\ x_7=v$$
 is a generalized f-variety, not locally Cohen-Macaulay, and $\dim V=4.$

Let $V_1$ be the variety defined by
$$x_1=s^4+t^4;\ x_2=s^2tu;\ x_3=s^3t;\ x_4=st^3;\ x_5=su^3;\
x_6=s^2t^2;\ x_7=v.$$
Let $K[V]$ and $K[V_1]$ be respectively  the coordinate rings of $V$
and  $V_1.$ It is immediate to check that  $V_1$ is a complete intersection,
 and therefore it is arithmetically Cohen-Macaulay,
   in fact $V_1$ is the Macaulayfication of $V$. 

\end{example}
 

\begin{thebibliography}{R-V}
 \bibitem{C} N.T. Cuong. Remarks on the Non-Cohen-Macaulay locus of Noetherian schemes.
{\em Proc. Amer. Math. Soc. 
} {\bf 
126-4} (1998),  1017--1022.  
\bibitem{Co} H. Cohn. A second course in number theory.
{\em New York and London: John Wiley and Sons, Inc.} XIII, 276 p. (1962).
 \bibitem{ES}
D.  Eisenbud and B.  Sturmfels.  Binomials ideals.  {\em Duke Math. J. 
} {\bf 
84-1} (1996),  1--45. 
 \bibitem{G} S.  Goto.  Approximatively 
Cohen-Macaulay Rings.  {\em Journal of Algebra} {\bf 76} (1982),  
214--225.  
\bibitem{GSW} S.  Goto,  N.  Suzuki and K.  Watanabe.  On 
affine semigroups.  {\em Japan Journal Math.} {\bf 2} (1976),  1--12.
\bibitem{H}
M.  Hochster.  Rings of invariants of tori, Cohen-Macaulay rings generated by monomials, and polytopes. 
{\em Ann. Math. }  {\bf 
96}, 318-337 (1972).    
\bibitem{K} T.  Kawasaki.  On Macaulayfication of Noetherian Schemes.  
{\em 
Trans.  AMS} {\bf 352-6} (2000,  2517--2552.
\bibitem{Mo1}
M.  Morales.  Equations des Vari\'et\'es Monomiales en codimension 
deux.  {\em Journal of Algebra} {\bf 175} (1995),  1082--1095.
\bibitem{PS}
I.  Peeva and B.  Sturmfels.  Syzygies of codimension 2 lattices 
ideals.  
{\em Maths Z.} {\bf 298-1} (1986),  145--167.
 \bibitem{TH} 
N.V.  Trung and L.T.  Hoa.  Affine Semigroups and Cohen-Macaulay 
rings 
generated by monomials.  {\em Trans. AMS} {\bf 298-1} (1986),  
145--167.
\bibitem{S1} P.  Schenzel.  Dualisierende komplexe in der lokalen 
Algebra und Buchsbaum-Ringe.  {\em Lectures Notes in Maths} {\bf 907} 
(1982),  
Springer-Verlag.
\bibitem{S2} P.  Schenzel.  On the use of local cohomology in Algebra 
and Geometry,  in Six lectures on commutative Algebra,  J.  Elias et 
als 
editors.  {\em Progress in Maths} (1998),  Birkhauser.
\bibitem{S3} P.  Schenzel. On birational Macaulayfications and Cohen-Macaulay canonical modules. 
{\em J. Algebra} {\bf 275, No.2}, 751--770 (2004).


\end{thebibliography}
 \end{document}